\documentclass[11pt]{article}
\usepackage{amsthm, amsmath, amssymb, amsfonts, url, booktabs, tikz, setspace, fancyhdr, bm}
\usepackage[hidelinks]{hyperref}
\usepackage{hyperref, enumerate}
\usepackage[shortlabels]{enumitem}
\usepackage[babel]{microtype}
\usepackage[english]{babel}
\usepackage[capitalise]{cleveref}
\usepackage{comment}
\usepackage{scalerel}
\usetikzlibrary{arrows.meta}

\usepackage[margin=20mm,bottom=13mm,footskip=7mm,top=13mm]{geometry}

\counterwithin{figure}{section}

\newtheorem{theorem}{Theorem}[section]
\newtheorem{proposition}[theorem]{Proposition}
\newtheorem{lemma}[theorem]{Lemma}
\newtheorem{observation}[theorem]{Observation}

\newtheorem{conjecture}[theorem]{Conjecture}
\newtheorem{claim}[theorem]{Claim}
\crefname{claim}{claim}{Claims}

\theoremstyle{definition}

\theoremstyle{remark}

\DeclareMathOperator{\dist}{dist}

\newcommand{\abs}[1]{\left\lvert#1\right\rvert}

\newcommand{\EE}{\mathbb{E}}
\newcommand{\PP}{\mathbb{P}}

\newcommand{\mA}{\mathcal{A}}
\newcommand{\mB}{\mathcal{B}}
\newcommand{\mC}{\mathcal{C}}

\newcommand{\D}{\Delta}

\newcommand{\e}{\varepsilon}
\newcommand{\g}{\gamma}
\newcommand{\eps}{\varepsilon}

\newcommand{\G}{\Gamma}

\def \E {\mathbb{E}}
\def \mE {\mathcal{E}}
\def \mL {\mathcal{L}}

\def \B {\mathcal{B}}

\def \R {\mathbb{R}}
\def \I {\boldsymbol{1}}

\def \bT {\bar{T}}
\def \hT {\hat{T}}
\def \hJ {\hat{J}}
\def \T {\tilde{T}}
\def \o {\H{o}}

\tikzstyle{p}+=[fill=black, circle, minimum width = 1pt, inner sep =
1pt]
\tikzstyle{w}+=[fill=white, draw, circle, minimum width = 1pt, inner sep =
1.5pt]

\newcounter{propcounter}

\begin{document}

\title{Approximate packing of independent transversals in locally sparse graphs}

\author{
Debsoumya Chakraborti\thanks{
Mathematics Institute, University of Warwick, Coventry, CV4 7AL, UK. Research supported by the European Research Council (ERC) under the European Union Horizon 2020 research and innovation programme (grant agreement No. 947978). 
E-mail: {\tt
debsoumya.chakraborti@warwick.ac.uk}.}
\and 
Tuan Tran\thanks{
School of Mathematical Sciences, University of Sciences and Technology of China, Hefei,
China. Research supported by National Key Research and Development Program of China 2023YFA1010201, and the Excellent Young Talents Program (Overseas) of the National Natural Science Foundation of China. 
E-mail: {\tt
trantuan@ustc.edu.cn}.}
}

\date{}

\maketitle

\begin{abstract}
Fix $\eps >0$ and consider a multipartite graph $G$ with maximum degree at most $(1-\eps)n$, parts $V_1,\ldots,V_k$ of the same size $n$, and where every vertex has at most $o(n)$ neighbors in any part $V_i$. Loh and Sudakov proved that any such $G$ has an independent transversal. 
They further conjectured that the vertex set of $G$ can be decomposed into pairwise disjoint independent transversals. In the present paper, we resolve this conjecture approximately by showing that $G$ contains $(1-\eps)n$ pairwise disjoint independent transversals.
As applications, we give approximate answers to questions of Yuster, and of Fischer, K\"uhn, and Osthus. 
\end{abstract}

\section{Introduction}\label{sec:intro}

\subsection{State of the art}

Given a multipartite graph $G$ with the vertex partition $V(G)=V_1\cup \cdots \cup V_k$, 
an \textit{independent transversal} of $G$ is an independent set in $G$ which contains exactly one vertex from each part $V_i$. The problem of finding sufficient conditions to ensure the existence of an independent transversal was asked over half a century ago by Bollob\'as, Erd\o s, and Szemer\'edi~\cite{BES75}. Since its inception, this question has generated much interest in the literature due to its relevance with other combinatorial notions such as linear arboricity, strong chromatic number, and list coloring, see for example \cite{ABZ07,A88,A92,H01,H04,J92,LS07,M03,R99,RS02,ST06,Y97CPC,Y97,Y21}.

Bollob\'as, Erd\o s, and Szemer\'edi~\cite{BES75} conjectured that if the parts $V_i$ have size at least twice the maximum degree of $G$, then there exists an independent transversal of~$G$. A simple application of the Lov\'asz Local Lemma achieves the following. Here, the constant $e$ refers to the base of the natural logarithm.
\begin{proposition} [Alon~\cite{A88}] \label{Alon}
Let $G$ be a multipartite graph with maximum degree $\D$, whose parts $V_1, \ldots, V_k$ all have size $|V_i| \ge 2e\D$. Then, $G$ has an independent transversal. 
\end{proposition}
Later, Haxell~\cite{H01} improved this factor $2e$ to $2$, resolving the original conjecture by an ingenious topological approach. This result is best possible by a construction with part sizes $2\D(G)-1$ due to \cite{J92,ST06,Y97}. Notably, this construction has high local degree. We refer the reader to \cite{H16} for a survey on such problems related to independent transversals.

We next briefly illustrate the importance of independent transversals in studying list coloring. 
Consider a graph $G$ on vertex set $V$ and a collection of lists of colors $\mL=\{L_v : v\in V\}$. A proper coloring of the vertex set is called \textit{$\mL$-coloring} if every vertex $v$ is colored by some color in~$L_v$. 
Reed~\cite{R99} conjectured the following about list coloring. 
If (i) for every vertex $v\in V$ and color $c\in L_v$, there are at most $\D$ neighbors $u$ of $v$ such that $c\in L_u$, and (ii) $\abs{L_v}\ge \D + 1$ for every $v\in V$, then $G$ admits an $\mL$-coloring.
If $\Delta$ is also the maximum degree of $G$, then a greedy coloring ensures the existence of an $\mL$-coloring. 
In general, Bohman and Holzman~\cite{BH02} disproved this conjecture. However, Reed's conjecture was asymptotically shown to be true by Reed and Sudakov~\cite{RS02} with lists of size at least $\D+o(\D)$.

The above list-coloring problem can be framed in the language of independent transversals by considering a $\abs{V}$-partite graph $\G$ with vertex partition $\bigcup_{v\in V} \{(v,c) : c\in L_v\}$, where two vertices $(v_1,c_1)$ and $(v_2,c_2)$ are connected by an edge in $\G$ if $v_1v_2$ is an edge in $G$ and $c_1=c_2$.
Notice that an independent transversal in $\G$ corresponds to an $\mL$-coloring of $G$. 
Moreover, the maximum degree of $\G$ is at most $\D$. 
Thus, Haxell's result~\cite{H01} directly implies Reed's conjecture when the lists have size at least $2\D$. 
Observe that every vertex in $\G$ has at most one neighbor in any part of $\G$.
This gives a special structure to the underlying graph $\G$, and Aharoni and Holzman (see~\cite{LS07}) speculated that these structures are enough to guarantee an independent transversal.
Generalizing the result of Reed and Sudakov~\cite{RS02}, Loh and Sudakov~\cite{LS07} established the following. 
For convenience, for a multipartite graph $G$ with the vertex partition $V_1\cup \cdots \cup V_k$, define the {\em local degree} of $G$ to be the maximum number of neighbors of a vertex in any part $V_i$.

\begin{theorem} [Loh--Sudakov~\cite{LS07}] \label{Loh}
For any $0 < \e < 1$, there exists $\g >0$ such that the following holds. Let $G$ be a multipartite graph with maximum degree at most $(1-\e)n$, parts $V_1,\ldots,V_k$ of size $|V_i|\ge n$, and local degree at most $\g n$. Then, it has an independent transversal.
\end{theorem}

The multiplicative factor $1-\e$ in the maximum degree condition is asymptotically optimal. For instance, let $G$ be the vertex-disjoint union of $n$ cliques~$K_{n+1}$, and let each part $V_i$ consist of precisely one vertex from each clique. The local degree is~$1$, the maximum degree is~$n$, and each part has size~$n$, but there is no independent transversal. 
For constructions achieving slightly better bounds, see \cite{haxell2023constructing}.
In more recent works, Glock and Sudakov \cite{GS22} and, independently, Kang and Kelly \cite{KK22} relaxed the maximum degree condition in \Cref{Loh} to an average degree condition.

\subsection{Main result}

Loh and Sudakov \cite{LS07} suggested that any graph $G$ that satisfies the assumptions of \Cref{Loh} has $n$ disjoint independent transversals. They proved the much weaker statement that maximum degree at most $n-o(n)$, local degree $o(n)$, and parts of size at least $2n$ are sufficient. In this paper, we prove an approximate version of their conjecture, thus strengthening \Cref{Loh}. Our proof combines arguments from \cite{GS22,LS07}, together with some additional ideas.

\begin{theorem}[Approximate packing of independent transversals] \label{thm}
For any $0<\e<1$, there exists $\g >0$ such that the following holds. Let $G$ be a multipartite graph with maximum degree at most $(1-\e)n$, parts $V_1,\ldots,V_k$ of size $|V_i| \ge n$, and local degree at most $\g n$. Then, $G$ contains $(1-\e)n$ pairwise disjoint independent transversals.
\end{theorem}

In \cite{GS22,KK22}, it is shown to be enough to have an average degree condition on every part instead of the maximum degree condition on the underlying graph along with the other hypotheses in \Cref{thm} to find a single independent transversal. However, with this weaker hypothesis, it is not possible to find an approximate decomposition. To see this, consider $0< \e< 1/12$ and $k=(1-\e)n$, and the graph $G$ obtained by taking a disjoint union of $kn/2$ isolated vertices and $n/4$ copies of $H$, where $H$ is the complete $k$-partite graph with two vertices in each part. $G$ has average degree at most $(1-\e)n$. There is a unique way (up to isomorphism) to represent $G$ as a $k$-partite graph with $n$ vertices in each part and with local degree at most two. Observe that any independent transversal must contain at most one vertex from each copy of $H$, thus must contain at least $k-n/4 = (3/4-\e)n$ isolated vertices. Thus, there can be at most $\frac{kn/2}{(3/4-\e)n}< (2/3+\e)n$ independent transversals in $G$. 

Now, we mention a couple of direct applications of our result.
Yuster~\cite{Y21} made the following conjecture. 
Let $2\le k\le n-1$ and let $G$ be a $k$-partite graph with vertex partition $V(G) = V_1 \cup \cdots \cup V_k$ such that $|V_i| = n$ for each $i \in [k]$. 
If every pair $(V_i,V_j)$ induces a perfect matching, then $G$ has a decomposition into independent transversals. An easy application of Hall's theorem proves it for $k\le n/2$, and Yuster proved it for $k\le n/(1.78)$. 
\Cref{thm} implies an approximate decomposition when $k\le n-o(n)$, giving an approximate answer to the conjecture. 

Next, we mention an application in packing list-colorings. As before, consider a graph $G$ on the vertex set $V$ and a collection of lists of colors $\mL=\{L_v : v\in V\}$. 
Cambie, Cames van Batenburg, Davies, and Kang~\cite{CCDK21} asked for the minimum number $D$ such that if (i) for every vertex $v\in V$ and color $c\in L_v$, there are at most $\D$ neighbors $u$ of $v$ such that $c\in L_u$, and (ii) $\abs{L_v}\ge D$ for every $v\in V$, then $G$ admits $D$ pairwise disjoint $\mL$-colorings. They conjectured that $D = \D + o(\D)$. Notice that, similar to before, this problem of packing list-colorings translates to the decomposition into independent transversals. Thus, if all the lists have size at least $\D + o(\D)$, then \Cref{thm} implies that $G$ admits $\D$ pairwise disjoint $\mL$-colorings, approximately resolving the conjecture.

\subsection{Multipartite Hajnal-Szemer\'edi theorem}
In this subsection, we demonstrate another application of \Cref{thm}. 
A perfect $K_k$-packing in a graph $G$ is a collection of vertex-disjoint $k$-cliques that covers all the vertices of $G$. Obviously, a necessary condition for a perfect $K_k$-packing in $G$ is that $k$ divides the number of vertices of $G$.
The fundamental result of Hajnal and Szemer\'edi~\cite{HS70,KK08} states that if $k$ divides $n$, then every $n$-vertex graph with minimum degree at least $(1-1/k)n$ contains a perfect $K_k$-packing. It is easy to see that the minimum degree condition is best possible by considering the balanced $n$-vertex complete $k$-partite graph. 
We now consider a multipartite analogue of this result of Hajnal and Szemer\'edi. Given a multipartite graph $G$, we define the {\em partite minimum degree} of $G$ to be the largest $d$ such that every vertex has at least $d$ neighbors in each part other than its own. 

\begin{conjecture}[Fischer~\cite{F99}, K\"uhn--Osthus~\cite{KO09}] \label{con:hajnal}
Let $k\ge 2$ and $G$ be a $k$-partite graph with parts $V_1,\ldots,V_k$ of the same size $n$. If the partite minimum degree of $G$ is at least $\left(1-\frac{1}{k}\right)n$, then $G$ has a perfect $K_k$-packing unless both $k$ and $n$ are odd and $G$ is isomorphic to a single exceptional graph. 
\end{conjecture}

This conjecture is straightforward for $k=2$. 
It was resolved for $k=3$ and sufficiently large $n$ by Magyar and Martin~\cite{MM02}, for $k=4$ and sufficiently large $n$ by Martin and Szemer\'edi~\cite{MS08}, and finally for every fixed $k$ and sufficiently large $n$ by Keevash and Mycroft~\cite{KM15JCTB}. 
Prior to this last result, an approximate version was established in~\cite{KM15,LM13}. 

\begin{theorem}[Keevash--Mycroft~\cite{KM15}, Lo--Markstr\"om~\cite{LM13}] \label{thm:hajnal}
For all $0<\delta <1$ and $k\ge 2$, there exists $n_0$ such that the following holds for every $n\ge n_0$. Let $G$ be a $k$-partite graph with parts $V_1,\ldots,V_k$ of the same size $n$. If the partite minimum degree of $G$ is at least $\left(1-\frac{1-\delta}{k}\right)n$, then $G$ has a perfect $K_k$-packing.
\end{theorem}

All previous progress to \Cref{con:hajnal} crucially assumes that $n$ is sufficiently large compared to~$k$. Without this assumption, we give an approximate answer to \Cref{con:hajnal}.

\begin{theorem} \label{thm:large r hajnal}
For all $0<\e,\delta <1$, there exists $n_0$ such that the following holds for every $n\ge n_0$. Let $k\ge 2$ and $G$ be a $k$-partite graph with parts $V_1,\ldots, V_k$ of the same size $n$. If $G$ has partite minimum degree at least $\left(1-\frac{1-\delta}{k}\right)n$, then it contains $(1-\e)n$ pairwise vertex-disjoint copies of $K_k$. 
\end{theorem}

\begin{proof}
We can (and do) assume $\e \le \delta$. 
Pick $\g=\g(\eps) >0$ according to \Cref{thm}, and for $k\ge 2$, pick $n_0(\delta, k)$ according to \Cref{thm:hajnal}. 
Let $k_0=\lceil \g^{-1}\rceil$ and $n_0=\max_{2\le k\le k_0} n_0(\delta, k)$.

When $2\le k\le k_0$, \Cref{thm:hajnal} implies that $G$ has $n$ pairwise vertex-disjoint copies of $K_k$. Now consider the case $k> k_0$. 
Define a multipartite graph $H$ with parts $V_1,\ldots,V_k$, where a pair of vertices $u,v$ forms an edge of $H$ if and only if they are from different parts and $uv \notin E(G)$.
Since the partite minimum degree of $G$ is at least $\left(1-\frac{1-\delta}{k}\right)n$, $H$ has local degree at most $\frac{1-\delta}{k}n\le \g n$ and maximum degree at most $k\cdot \frac{1-\delta}{k}n\le (1-\e)n$. Thus, by \Cref{thm}, there are $(1-\e)n$ pairwise disjoint independent transversals in $H$, which translates to $(1-\e)n$ pairwise vertex-disjoint copies of $K_k$ in $G$.
\end{proof}

\subsection{Proof strategy}
In this section, we give a high-level overview of the proof of \Cref{thm}. To do this, we build upon the proof of \cite{GS22,LS07}. We begin with sketching their proof of \Cref{Loh}. 

Suppose we are given a $k$-partite graph $G$ with the vertex partition $V_1 \cup \cdots \cup V_k$, where each part $V_i$ contains $n$ vertices and the maximum degree of $G$ is at most $(1-\e)n$. Notice that if the maximum degree of $G$ is at most $n/(2e)$, then the existence of an independent transversal follows from \Cref{Alon}. Loh and Sudakov \cite{LS07} utilized this together with a randomized nibble algorithm to build a transversal in several iterations, where each iteration adds a few vertices to the partial transversal made so far. Crucially, their algorithm ensures that the maximum degree decreases faster than the size of the remaining vertices in each part so that after some iterations, the ratio between the size of each part and the maximum degree of the graph becomes at least $2e$. 

We next show a heuristic of why the first iteration of the above algorithm works. We activate each part with some small probability $p$ and then choose a vertex uniformly at random from each activated part. Denote this set of all chosen vertices by $T$; we add these vertices to the transversal unless they have a neighbor in $T$. To avoid conflicts in later iterations, we delete all the neighbors of the vertices in $T$. Thus, for any vertex, by a simple union bound, the probability that one of its neighbors is chosen in $T$ is at most $(1-\e) n\cdot p/n = (1-\e)p$. Thus, we expect to lose at most a $(1-\e)p$-fraction of vertices in each part. Moreover, we discard all the parts from which a vertex has already been added to the transversal. Clearly, each part $V_i$ is deleted when $T$ contains a vertex $v\in V_i$, but there are no neighbors of $v$ in $T$. This happens with a probability of at least $p(1-(1-\e)p)\ge p-p^2$ that results in at least a $(p-p^2)$-fraction decay in the degrees of the vertices. This shows that the maximum degree is reduced by a larger factor than the part sizes for sufficiently small $p$. The accumulation of these factors through several iterations gives us the desired factor of $2e$ for the application of \Cref{Alon}.    

In view of this method, to prove \Cref{thm}, a naive approach is to repeatedly apply the same argument to find an independent transversal and then remove it from the graph. However, such a naive approach does not allow us to control the maximum degree along this process. For example, after removing a single transversal from $G$, the maximum degree may not be reduced. To see this, let $k\ge 2(1-\e)n$ and consider any $k$-partite graph $G$ with maximum degree $(1-\e)n$, parts of size at least $n$ such that $G$ contains a copy $S$ of $K_{(1-\e)n,(1-\e)n}$, where every pair of vertices in $S$ lie in different parts of $G$. Firstly, it is unlikely for the above randomized algorithm to add one of the entire parts of the bipartition of $S$ to the transversal. Secondly, the algorithm keeps one of the parts completely untouched as the transversal must form an independent set. Thus, the maximum degree remains the same after removing a transversal from such $G$. Since the above nibble scheme heavily relies on the fact that the maximum degree is lower than the part sizes, 
\begin{itemize}[leftmargin=*]
\item[$\star$] we need to be able to remove transversals while ensuring that the maximum degree of the remaining graph remains lower than the part sizes. 
\end{itemize}

To prove \Cref{thm}, instead of sequentially building the transversals, another extreme approach is to build all transversals in parallel. However, the way a single transversal is completed in \cite{LS07} using \Cref{Alon} does not seem to extend easily to build $(1-\eps)n$ pairwise disjoint independent transversals at once. Thus, to facilitate the application of \Cref{Alon} while maintaining the property in $\star$, we run a randomized algorithm for several rounds and build a few disjoint transversals in a single round. The main idea of this paper is to add another layer of nibble algorithm to the original algorithm of \cite{LS07}. This allows us to control the maximum degree in future rounds as desired. This extends the standard nibble method to a two-dimensional form, where a single round consists of constructing some disjoint transversals, and on the other hand, each transversal of a single round is built in several iterations, where each iteration adds a few vertices to all the transversals. This way of using the nibble method may have other potential applications in graph decomposition problems.

It is more convenient to use concentration inequalities to analyze our random process if the local degree is a much lower order term in $n$ in \Cref{thm}. Fortunately, we have some standard techniques to reduce our main result to the case where the local degree is just a constant.
Instead of dealing with \Cref{thm} directly, we will first use our randomized algorithm to prove the following theorem and later deduce \Cref{thm} from it. 

\begin{theorem} \label{thm:main}
For any $0<\e<1$ and $C>0$, there exists $n_0$ such that the following holds for all $n\ge n_0$. Let $G$ be a multipartite graph with maximum degree at most $(1-\e)n$, parts $V_1,\ldots, V_k$ of size $|V_i| \ge n$, and local degree at most $C$. Then, it contains $(1-\e)n$ pairwise disjoint independent transversals.
\end{theorem}

\medskip
\noindent{\bf Organization.} The remainder of this paper is organized as follows. In the next section, we mention several probabilistic tools that will be useful to us. 
The subsequent three sections are devoted to proving \Cref{thm:main}. For that, we start by giving a framework in \Cref{sec:frame} and introducing a few parameters that we will keep track of in our randomized process. Then, we describe the randomized algorithm in \Cref{sec:algorithm} and analyze this algorithm in \Cref{sec:analysis} to finish the proof of \Cref{thm:main}. 
In \Cref{sec:reduction}, we show how to deduce \Cref{thm} from \Cref{thm:main}. 
Finally, we finish with a couple of concluding remarks in \Cref{sec:remarks}.

\medskip
\noindent{\bf Notation.} For positive integers $n$, we use $[n]$ to denote the set $\{1,2,\ldots,n\}$. For nonnegative reals $a,b,c$, we write $a= b\pm c$ to mean that $b-c \le a \le b+c$. 
Unless stated otherwise, all asymptotics are as $n\rightarrow \infty$. We say that an event $\mA$ occurs \textit{with very high probability} (\textit{w.v.h.p.} for short) if $\PP[\mA]\ge 1-n^{-\omega(1)}$. It is convenient to work with this notion due to the fact that if we have a collection of $n^{O(1)}$ events where each event holds w.v.h.p., then their intersection also holds w.v.h.p. 
For a simple presentation, we omit the rounding signs throughout the paper. 

\medskip
We use standard graph theoretic notations. 
Consider a graph $G$. We denote the vertex set of $G$ by $V(G)$, the edge set by $E(G)$, and the maximum degree by $\D(G)$. 
For $V\subseteq V(G)$, we denote by $G[V]$ the subgraph of $G$ induced by $V$. 
For a vertex $v\in V(G)$, we denote the degree of $v$ by $d_G(v)$ and the neighborhood of $v$ by $N_G(v)$. 
For $v\in V(G)$ and $V\subseteq V(G)$, we denote by $d_G(v,V)$ the size of the set $V\cap N_G(v)$.
For $V\subseteq V(G)$, we denote $N_G(V) = \bigcup_{v\in V} N_G(v)$. 
For $u,v\in V(G)$, the distance between $u$ and $v$, denoted by $\dist_G(u,v)$, is the minimum number of edges in a path in $G$ between them. 
When the underlying graph $G$ is clear from context, we may omit the subscript from these notations.

\section{Probabilistic tools}
To analyze our random process, we will use a few concentration inequalities.

\begin{lemma}[{Chernoff's bound, see~\cite{MU05}}] \label{chernoff}
Let $X$ be the sum of independent Bernoulli random variables (possibly with distinct expectations). Then, for all $0\le s \le \E(X)$,
\[
\PP\left[\abs{X-\E(X)}\ge s\right]\le 2 \exp\left(-\frac{s^2}{3\E(X)}\right).
\]
Moreover, for $s\ge 7\E(X)$, we have $\PP[X\ge s] \le e^{-s}$. 
\end{lemma}

\begin{lemma}[Hoeffding's inequality, see~\cite{H63}]\label{hoeffding} 
Let $X=X_1+\cdots +X_n$, where the $X_i$'s are independent random variables such that $0\leq X_i\leq L$ for all $i$. Then, for every $s>0$, we have
\[
\PP\left[\abs{X-\E(X)}\ge s\right]\le 2\exp\left(-\frac{s^2}{2L^2n}\right).
\]
\end{lemma}

To state the next inequality, we need a couple of definitions. 
Given a product probability space $\Omega = \prod_{i=1}^n \Omega_i$, we say the following. 
\begin{itemize}[leftmargin=*]
\item A random variable $X\colon \Omega \rightarrow \R$ is $L$-Lipschitz, if for every $\omega \in \Omega$, changing $\omega$ in any single coordinate affects the value of $X(\omega)$ by at most $L$.
\item A non-negative random variable $X\colon \Omega \rightarrow \R_{+}$ is $r$-certifiable, if for every $\omega \in \Omega$ and $s\ge 0$ with $X(\omega) \ge s$, there exists a set $I \subseteq \{1,\ldots,n\}$ of size at most $rs$ such that $X(\omega') \ge s$ for every $\omega' \in \Omega$ that agrees with $\omega$ on the coordinates indexed by $I$.
\end{itemize}

\begin{lemma} [Talagrand's inequality, see~\cite{MR02}] \label{talagrand}
Suppose that $X$ is an $L$-Lipschitz and $r$-certifiable random variable. Then, for every $0\le s\le \E(X)$, we have
\[
\PP\left[\abs{X - \E(X)} > s + 60L\sqrt{r\E(X)}\right] \le 4\exp\left(-\frac{s^2}{8L^2r\E(X)}\right).
\]
\end{lemma}

Since the number of parts $k$ in \Cref{thm} can be arbitrarily large compared to $n$, simple union bounds do not work while analyzing our random process. 
In such situations, we will use the Lov\'asz Local Lemma which allows us to prove that with a positive probability, none of the ``bad'' events happen, provided that the events have limited dependency. Given events $\B_1,\ldots,\B_n$, a graph $\G$ on $[n]$ is a {\em dependency graph} for these events if each event $\B_i$ is mutually independent of all other events except those indexed by $N_{\G}(i)\cup \{i\}$.

\begin{lemma} [Lov\'asz Local Lemma, see \cite{AS16}] \label{lll}
Let $\B_1, \ldots, \B_n$ be events in a probability space with dependency graph $\G$. If $\PP[\B_i] \le p$ for all $i\in [n]$ and $ep(\D(\G)+1)\le 1$, then $\PP[\cap \overline{\B_i}] > 0$.
\end{lemma}

\section{Frame of the proof} \label{sec:frame}
This section and the next two sections contain the proof of \Cref{thm:main}. In this section, we introduce many parameters that we would like to keep track of throughout our random process. We can (and do) assume that $\e$ is a sufficiently small positive constant and $n$ is sufficiently large with respect to $\e$ to support our arguments. Since deleting vertices from $G$ does not increase the maximum degree or the local degree, without loss of generality, we assume $|V_i| = n$ for each $i \in [k]$. Set $p=1/\log^3 n$.

In every round, we plan to make some number of independent transversals. 
If after some rounds the part sizes are $n'$, then in the next round, we simultaneously construct $pn'$ transversals. 
Throughout this paper, we use $r$ to track the number of rounds and use $t$ to track the number of iterations in a single round.
Set $\delta=\eps^5$ and $r^*=t^*=\frac{30}{\e p}$. We will use $r^*$ to denote the total number of rounds, and $t^*$ to denote the total number of iterations in each round. We define the following sequences to keep track of part sizes and maximum degrees throughout the randomized algorithm.

\begin{itemize}[leftmargin=*]
\item $S_r(0)=(1-p)^{r-1}n$, $D_r(0)=(1-\e) \left(1-p+\e^3 p\right)^{r-1}n$, and $p_r=\frac{D_r(0)p}{S_r(0)}$.
\item $S_r^{-}(t)=(1-p_r-p^2)^t S_r(0)$ and $S_r^{+}(t)=(1-p_r+p^2)^t S_r(0)$. In particular, $S_r^{-}(0)=S_r^{+}(0)=S_r(0)$.
\item $D_r(t) = \left(1 - p+ \e p/2\right)^t D_r(0)$. 
\end{itemize}
We now briefly describe how the above parameters will appear in our randomized algorithm. After removing the independent transversals constructed until the start of the $r$-th round, in the remainder of the graph, the part sizes will be around $S_r(0)$, and the maximum degree will be at most $D_r(0)$. One can check that this is readily satisfied when $r=1$ by the assumptions in \Cref{thm:main}. The quantities $S_r^{-}(t)$, $S_r^{+}(t)$, and $D_r(t)$ will help us control the part sizes and maximum degree after the $t$-th iteration during the $r$-th round.
We will use the following relations in our proof.

\begin{observation}\label{obs} The following hold for all $r\in [r^*]$ and $0\le t \le t^*$.
\begin{itemize}
    \item[(i)] $(1-3p)^t=\Omega(1)$; 
    \item[(ii)] $(1-\eps)p\le p_r\le (1-2\eps/3)p$;
    \item[(iii)] $\frac{D_r(t)}{S^{-}_r(t)} \le \frac{D_r(0)}{S_r(0)} \le 1$;
    \item[(iv)] $n\ge S^{-}_r(t)\ge D_r(t)=\Omega(n)$;
    \item[(v)] $\frac{D_r(t^*)}{S^{-}_r(t^*)-pS_r(0)} < \frac{1}{2e}$.
\end{itemize}
\end{observation}
\begin{proof}
(i) Using the inequality $1-x \ge e^{-2x}$ for $x \in [0,\tfrac12]$ yields 
$(1-3p)^t \ge e^{-6pt} \ge e^{-6pt^*}=e^{-180/\eps}=\Omega(1).$

(ii) We see that $\frac{D_r(0)}{S_r(0)}=(1-\eps)\left(\frac{1-p+\eps^3p}{1-p}\right)^{r-1} \ge 1-\eps$, and that
\[
\frac{D_r(0)}{S_r(0)}=(1-\eps)\left(\frac{1-p+\eps^3p}{1-p}\right)^{r-1} \le (1-\eps)(1+2\eps^3 p)^{r-1} \le (1-\eps) e^{2\eps^3 p (r-1)}\le (1-\eps)e^{60\eps^2} \le 1-2\eps/3.
\]
Combining with the fact that $p_r=\frac{D_r(0)p}{S_r(0)}$, we get $(1-\eps)p\le p_r\le (1-2\eps/3)p$.

(iii) As $p_r \le (1-2\eps/3)p$, we find $1-p_r-p^2 \ge 1-p+2\eps p/3-p^2 \ge 1-p+\eps p/2$. Thus, 
\[
\frac{D_r(t)}{S^{-}_r(t)}=\frac{D_r(0)}{S_r(0)}\left(\frac{1-p+\eps p/2}{1-p_r-p^2}\right)^{t} \le \frac{D_r(0)}{S_r(0)} \overset{(ii)}{\le} 1.
\]

(iv) Since $0\le p_r \le p=o(1)$, we have $S^{-}_r(t)=(1-p_r-p^2)^t(1-p)^{r-1}n \le n$. 
In addition, since both $1-p+\eps^3 p$ and $1-p+\eps p/2$ are greater than $1-p$, we get 
\[
D_r(t)\ge (1-\eps) (1-p)^{r+t-1}n \ge \tfrac12 e^{-2p(r+t-1)}n \ge \tfrac12 e^{-120/\eps}n=\Omega(n).
\]
Therefore, combining with (iii), we conclude $n\ge S^{-}_r(t)\ge D_r(t)=\Omega(n)$.

(v) We have
\begin{equation*}
\frac{D_r(t^*)}{S^{-}_r(t^*)} =\frac{D_r(0)}{S_r(0)}\left(\frac{1-p+\eps p/2}{1-p_r-p^2}\right)^{t^*} \overset{(ii),(iii)}{\le} 
\left(\frac{1-p+\eps p/2}{1-p+2\eps p/3-p^2}\right)^{t^*}\le \left(1-\eps p/7\right)^{t^*}\le e^{-\eps p t^*/7}< \frac{1}{4e}.
\end{equation*}
Hence, $\frac{D_r(t^*)}{S^{-}_r(t^*)-pS_r(0)} \overset{(i),(iii)}{\le} \frac{D_r(t^*)}{(1/2)S^{-}_r(t^*)} < \frac{1}{2e}$.
\end{proof}

During the $r$-th round, we will construct $pS_r(0)$ independent transversals through several iterations. As we proceed with these iterations, these independent transversals will only be partially built. We use $\ell$ to denote such partial transversals, and the set of these partial transversals will be denoted by~$[p S_r(0)]$.
We will inductively ensure the existence of the following during the execution of our algorithm. For every ``round'' $r\in [r^*+1]$, we guarantee the following. 
\stepcounter{propcounter}
\begin{enumerate}[leftmargin=*,label = {\bfseries \emph{\Alph{propcounter}\arabic{enumi}}}]
    \item 
    There are $pS_1(0)+\dots + pS_{r-1}(0)$
    pairwise disjoint independent transversals in $G$. Let $L_r$ denote the set of all vertices used in these transversals. \label{ls}
    \item The maximum degree of the graph $G_r := G\left[\bigcup_{i\in [k]}V_i\setminus L_r\right]$ is at most $D_r(0)$. 
    \label{vij0}
\end{enumerate}
For convenience, for every $i\in [k]$ and ``partial independent transversal'' $\ell\in [p S_r(0)]$, define $V_i^{\ell}(0)= V_i\setminus L_r$. Observe that $\abs{V_i^{\ell}(0)}= S_r(0)$ for all $i\in [k]$ and $\ell\in [p S_r(0)]$.
Furthermore, during every ``round'' $r\in [r^*]$, for each ``partial independent transversal'' $\ell\in [p S_r(0)]$ and ``iteration'' $t\in [t^*]$, we make sure the existence of the objects introduced in \ref{isjt}--\ref{vijt} such that the properties in \ref{lower bound vijt}--\ref{neighbor in remaining parts} hold.

\stepcounter{propcounter}
\begin{enumerate}[leftmargin=*,label = {\bfseries \emph{\Alph{propcounter}\arabic{enumi}}}]
    \item There is $I^{\ell}_r(t)\subseteq [k]$. (As we will see below, for the $\ell$-th transversal of the $r$-th round, $I^{\ell}_r(t)$ is the set of indices of partition classes still need to be visited after the $t$-th iteration.) \label{isjt} 
    \item There is an independent transversal $T^{\ell}(t)$ of $G\left[\bigcup_{i\in [k]\setminus I^{\ell}_r(t)}V_i\right]$. Moreover, for every fixed $t\in [t^*]$, the sets in $\{T^{\ell}(t) : \ell\in [p S_r(0)]\}$ are pairwise disjoint. Define $T(t)=\bigcup_{\ell\in [pS_r(0)]} T^{\ell}(t)$. \label{tjt}
    \item There is a subset $V_i^{\ell}(t)\subseteq V_i^{\ell}(0)\setminus N(T^{\ell}(t))$ for each $i\in I^{\ell}_r(t)$. Define $V^{\ell}(t)=\bigcup_{i\in I_r^{\ell}(t)} V_i^{\ell}(t)$ and $G^{\ell}(t) = G[V^{\ell}(t)]$. (The set $V_i^{\ell}(t)$ should be thought of as the remaining vertex set in the $i$-th part for the $\ell$-th transversal after the $t$-th iteration. Note that since the set $V^{\ell}(t)$ do not contain any neighbors of $T^{\ell}(t)$, it is sufficient to find an independent transversal of $G^{\ell}(t)$ in order to complete $T^{\ell}(t)$ to an independent transversal of $G$.) \label{vijt} 
\end{enumerate}
\stepcounter{propcounter}
\begin{enumerate}[leftmargin=*,label = {\bfseries \emph{\Alph{propcounter}\arabic{enumi}}}]
    \item For every $\ell \in [p S_r(0)]$ and $i\in I^{\ell}_r(t)$, we have $S_r^{-}(t)\le \abs{V_i^{\ell}(t)}\le S_r^{+}(t)$. \label{lower bound vijt} 
    \item Each $v\in V(G_r)$ has at most $D_r(t)$ neighbors in $V^{\ell}(t)$. In particular, $G^{\ell}(t)$ has maximum degree at most $D_r(t)$. \label{upper bound degree} 
    \item For every $v\in V(G_r)$ with $d_{G_r}(v)\ge \log^{15} n$, 
    \begin{equation*}\label{neighbor in transversals}
    \abs{N_{G_r}(v) \cap T(t)} \ge  (1-3\delta)p^2 \sum_{0\le j< t} (1-(1+2\delta)p)^j d_{G_r}(v).
    \end{equation*}
    \item For every $v\in V(G_r)$, 
    \begin{equation*}\label{neighbor in remaining parts}
    \sum_{\ell\in [p S_r(0)]}\abs{N_{G_r}(v) \cap V^{\ell}(t)} \ge pS_r(0)\cdot (1- p-p_r-\delta p)^t d_{G_r}(v).
    \end{equation*}
\end{enumerate}

Considering $pS_r(0)$ many transversals simultaneously plays a crucial role in establishing the property~\ref{neighbor in remaining parts} (as also remarked after \Cref{claim:neighbor in remaining parts}), which itself is a key step in our proof of \Cref{thm:main}. As the readers will see later, the property~\ref{neighbor in remaining parts} serves as an intermediate step for reaching the property \ref{neighbor in transversals} which is one of our main means to control the maximum degree as the rounds progress. Roughly speaking, it asserts that each vertex has many neighbors in the union of already built transversals, and so after removing these transversals from $G_r$, each vertex loses the right proportion of its neighbors.

We start by showing that for every $r\in [r^*]$, during the $r$-th round, after the $t^*$-th iteration, we can extend the partial independent transversals $T^{\ell}(t^*)$ to pairwise disjoint full independent transversals one by one by using \Cref{Alon}. Indeed, suppose for some $\ell\in [pS_r(0)]$, we have constructed $\ell-1$ pairwise disjoint independent transversals by extending $T^{1}(t^*),\dots,T^{\ell-1}(t^*)$. Then, the graph obtained from $G^{\ell+1}(t^*)$ after removing the vertices in these $\ell-1$ transversals, has part sizes at least $S^{-}_r(t^*) - pS_r(0)$ and has maximum degree at most $D_r(t^*)$ (by \ref{lower bound vijt} and \ref{upper bound degree}). Thus, by \Cref{Alon} and \Cref{obs}~(v), we can extend $T^{\ell}(t^*)$ to an independent transversal that is disjoint from the already built $\ell-1$ transversals.
Thus, assuming that the objects in \ref{isjt}--\ref{vijt} satisfy \ref{lower bound vijt}--\ref{neighbor in remaining parts} after the $t^*$-th iteration of the $r$-th round, we can move to the round $r+1$ and \ref{ls} is already satisfied. Moreover, \ref{neighbor in transversals} implies that \ref{vij0} holds for $r+1$ as follows. To see that, consider any vertex $v\in V(G_r)$. If $d_{G_r}(v) \le \log^{15} n$, then $d_{G_r}(v) \le \log^{15} n \le D_{r+1}(0)$. Otherwise, by \ref{neighbor in transversals}, we see that the degree of $v$ in $G_{r+1}$ will be at most 
\begin{align*}
d_{G_r}(v) - \abs{N_{G_r}(v)\cap T(t^*)} &\le D_r(0) - (1-3\delta)p^2 \sum_{0\le j< t^*} (1-(1+2\delta)p)^j D_r(0) \displaybreak[0] \\
&=D_r(0)\Big(1-\frac{p(1-3\delta)}{1+2\delta}\left[1-(1-(1+2\delta)p)^{t^*}\right]\Big) \displaybreak[0] \\
& \le D_r(0)\Big(1-\frac{p(1-3\delta)}{1+2\delta}\left[1-e^{-(1+2\delta)pt^*}\right]\Big) \displaybreak[0] \\
&\le D_r(0)\left(1-p+\eps^3 p\right)=D_{r+1}(0),
\end{align*}
where the last inequality holds since $\frac{(1-3\delta)(1-e^{-(1+2\delta)pt^*})}{1+2\delta} =\frac{(1-3\delta)(1-e^{-30(1+2\delta)/\eps})}{1+2\delta} \ge \frac{(1-3\delta)(1-\eps^4)}{1+2\delta}\ge 1-\eps^3$ assuming that $\eps>0$ is sufficiently small (recall that $\delta=\eps^5$).

Note that if we are able to ensure \ref{ls} at the end of the $r^*$-th round, then the number of pairwise disjoint independent transversals in $G$ is given by
\[
p\sum_{i\in [r^*]} S_i(0) = pn\sum_{i\in [r^*]} (1-p)^{i-1} = n\left(1-(1-p)^{r^*}\right)\ge n\left(1-e^{-pr^*}\right) = n\left(1-e^{-30/\eps}\right)\ge (1-\eps) n,
\]
where the last inequality uses the fact that $\eps$ is sufficiently small. 
Thus, to prove \Cref{thm:main}, the only thing remaining is to verify \ref{isjt}--\ref{vijt} and \ref{lower bound vijt}--\ref{neighbor in remaining parts} for the $(t+1)$-st iteration assuming the same hold for some $0\le t < t^*$. Note that these properties trivially hold for $t=0$ by setting $I_r^{\ell}(0)=[k]$ and $T^{\ell}(0)=\emptyset$, and by using~\ref{ls} and~\ref{vij0} .

\section{Randomized algorithm} \label{sec:algorithm}
Suppose we are at the $r$-th round for some $r\in [r^*]$. Suppose for some $0\le t< t^*$, we completed the $t$-th iteration with the objects in \ref{isjt}--\ref{vijt} satisfying \ref{lower bound vijt}--\ref{neighbor in remaining parts}. We next describe the $(t+1)$-st iteration during the $r$-th round (however, for convenience, we often drop $r,t$ from the notations). For every $i\in I_r^{\ell}(t)$ and $\ell\in [p S_r(0)]$, suppose $V_i^{\ell}(t)$ is the remaining vertex set for the $\ell$-th transversal after the $t$-th iteration as introduced in \ref{vijt}. 
\begin{enumerate}[leftmargin=*]
    \item For every $\ell\in [p S_r(0)]$ and $i\in I^{\ell}_r(t)$, uniformly at random select a vertex $v_i^{\ell}$ from $V_i^{\ell}(t)$. 
    For $\ell\in [p S_r(0)]$, let $\T^{\ell} =\{v_i^{\ell}: i\in I^{\ell}_r(t)\}$ be the set of all selected vertices. Note that $\T^{\ell}$ might induce edges.
    \item For every $\ell\in [p S_r(0)]$ and $i\in I_r^{\ell}(t)$, activate each part $V_i^{\ell}(t)$ independently with probability~$p$. Let $\hJ^{\ell}$ be the set of all indices corresponding to the activated parts for the $\ell$-th transversal, that is, $\hJ^{\ell}=\{i\in I_r^{\ell}(t): V_i^{\ell}(t) \; \text{is activated}\}$.
    Let $\hT^{\ell} =\{v_i^{\ell}: i\in \hJ^{\ell}\}$ be the set of all selected vertices whose parts are activated. Denote by $\hT$ the multi-set (a vertex can be picked more than once) by concatenating~$\hT^{\ell}$'s. 
    \item For each $\ell\in [p S_r(0)]$, let $J^{\ell}$ be the set of all indices $i\in \hJ^{\ell}$ such that $v_i^{\ell}$ appears exactly once in $\hT$ and it is not adjacent to any other vertex of $\hT^{\ell}$. Denote the collection of these vertices by $\bT^{\ell}$, that is, $\bT^{\ell}=\{v_i^{\ell}: i\in J^{\ell}\}$. Observe that each $\bT^{\ell}$ forms an independent set, and the $\bT^{\ell}$'s are pairwise disjoint. We thus, for each $\ell$, add the vertices in $\bT^{\ell}$ to the $\ell$-th partial independent transversal in the $r$-th round, that is, we let $T^{\ell}(t+1)= T^{\ell}(t)\cup \bT^{\ell}$. 
    \item For each $i\in J^{\ell}$, delete the entire part $V_i^{\ell}(t)$ corresponding to $i$ from consideration for the $\ell$-th transversal in the $r$-th round. In other words, set $I_r^{\ell}(t+1)=I_r^{\ell}(t)\setminus J^{\ell}$. 
    This, together with the last step, defines the objects in \ref{isjt} and \ref{tjt} for the $(t+1)$-st iteration.
    \item We next aim to define the sets $V_i^{\ell}(t+1)$ for $i\in I_r^{\ell}(t+1)$ and $\ell\in [pS_r(0)]$. As per \ref{vijt}, we need to ensure $V_i^{\ell}(t+1)\subseteq V_i^{\ell}(0)\setminus N(T^{\ell}(t+1))$.
    Thus, to construct $V_i^{\ell}(t+1)$ from $V_i^{\ell}(t)$, we will delete all neighbors of the vertices in $\bT^{\ell}$. For the convenience of the analysis, we potentially delete more vertices than that. 
    Observe that for every vertex $v\in V^{\ell}(t)$, each of its remaining neighbors is included in $\hT^{\ell}$ with probability at most $\frac{p}{S^{-}_r(t)}$.
    Thus, by a simple union bound, we have 
    \[
    \PP[N_{G^{\ell}(t)}(v)\cap \hT^{\ell} \neq \emptyset]\le d_{G^{\ell}(t)}(v)\cdot \frac{p}{S^{-}_r(t)}\overset{\text{\ref{upper bound degree}}}{\le} \frac{D_r(t)p}{S^{-}_r(t)} \overset{\text{Obs. \ref{obs}(iii)}}{\le} \frac{D_r(0)p}{S_r(0)}= p_r.
    \]
    For every $i\in I_r^{\ell}(t+1)$, while constructing $V_i^{\ell}(t+1)$ from $V_i^{\ell}(t)$, in parallel with deleting all vertices in $N_{G^{\ell}(t)}(\hT^{\ell})$, we artificially delete all the vertices $v\in V_i^{\ell}(t)$ according to a Bernoulli random variable $B^{\ell}_v$, independent of all other variables, so that 
    \begin{align} \label{eq:same probability}
       \PP[N_{G^{\ell}(t)}(v)\cap \hT^{\ell} \neq \emptyset \;\; \text{or} \;\; B^{\ell}_v = 1] = p_r. 
    \end{align}
    This ensures that in the $(t+1)$-st iteration, every vertex is deleted with the same probability. 
    \item For $\ell\in [pS_r(0)]$ and $i\in I_r^{\ell}(t+1)$, let $V_i^{\ell}(t+1)$ denote the set of all vertices $v$ in $V_i^{\ell}(t)$ not satisfying the event in \eqref{eq:same probability} (i.e. $N_{G^{\ell}(t)}(v)\cap \hT^{\ell} = \emptyset$ and $B^{\ell}_v = 0$).
    Since $\bT^{\ell}\subseteq \hT^{\ell}$, we have $V_i^{\ell}(t+1)\subseteq V_i^{\ell}(0)\setminus N(T^{\ell}(t+1))$ for every $i\in I_r^{\ell}(t+1)$. For every $\ell\in [pS_r(0)]$, the sets $V_i^{\ell}(t+1)$ along with the definitions of \ref{vijt} introduces $V^{\ell}(t+1)$ and $G^{\ell}(t+1)$. 
\end{enumerate}
In this section, we already defined the objects in \ref{isjt}--\ref{vijt} for the $(t+1)$-st iteration of the $r$-th round. In the next section, we show that the properties \ref{lower bound vijt}--\ref{neighbor in remaining parts} also hold simultaneously with positive probability for the $(t+1)$-st iteration. This will finish the proof of \Cref{thm:main}.

\section{Analysis of the algorithm} \label{sec:analysis}
Fix $r\in [r^*]$ and $0\le t< t^*$. As done at the beginning of the last section, we assume that the objects in \ref{isjt}--\ref{vijt} satisfy the properties \ref{lower bound vijt}--\ref{neighbor in remaining parts} for these fixed $r$ and $t$. Our goal is to show that there is a choice of the objects defined in the randomized algorithm for the $(t+1)$-st iteration of the $r$-th round so that the properties remain true. To achieve this, we first show that the events for individual parts or vertices corresponding to \ref{lower bound vijt}--\ref{neighbor in remaining parts} hold with very high probability, and then by the local lemma we show that they simultaneously occur with positive probability.

\subsection{Size of remaining parts}
In the following, we show that for fixed $\ell$ and $i$, the event \ref{lower bound vijt} holds w.v.h.p. 

\begin{claim}\label{claim:B4}
For any $\ell \in [p S_r(0)]$ and $i\in I^{\ell}_r(t+1)$, w.v.h.p. we have 
\[
\abs{V_i^{\ell}(t+1)} = (1-p_r)\abs{V_i^{\ell}(t)} \pm \log n \sqrt{p_r\abs{V_i^{\ell}(t)}}
\]
\[
\text{and} \;\;\;\;\;\;\;\; S_r^{-}(t+1)\le \abs{V_i^{\ell}(t+1)}\le S_r^{+}(t+1).
\]
\end{claim}

\begin{proof}
Let $R=\abs{V_i^{\ell}(t)\setminus V_i^{\ell}(t+1)}$ be the number of vertices removed from $V_i^{\ell}(t)$. Since the local degree of $G$ is bounded by $C$, changing the outcome of the selected vertex in a single part in Step~1 can affect $R$ by at most $C$, and changing the outcome of activating a single part in Step~2 can affect $R$ by at most~$C$. Moreover, for any vertex $v$ removed, there must be either a neighbor of $v$ in $\hT^{\ell}$, or we have $B^{\ell}_v=0$. Therefore, $R$ is $C$-Lipschitz and $1$-certifiable. Note that $\E(R) \overset{\eqref{eq:same probability}}{=} p_r\abs{V_i^{\ell}(t)}$. Thus, applying Talagrand's inequality with $s=\frac{\log n}{2} \sqrt{p_r\abs{V_i^{\ell}(t)}}$ (it can be easily checked that $s\le \E(R)$ from the assumption that $S_r^{-}(t)\le \abs{V_i^{\ell}(t)}$ and \Cref{obs} (ii), (iv)),
we obtain 
\[
\PP\left[|R-\EE(R)|>2s\right]\le \exp\left(-\Omega(\log^2 n)\right),
\]
since $60C\sqrt{\EE(R)} = 60C\sqrt{p_r\abs{V^{\ell}_i(t)}} \le s$ and $s^2/\EE(R) = (\log^2 n)/4$. This proves the first part of the claim as $\abs{V_i^{\ell}(t+1)}-(1-p_r)\abs{V_i^{\ell}(t)}=\EE(R)-R$.

For the second part of the claim, note that by assumption, we have $S_r^{-}(t)\le \abs{V_i^{\ell}(t)}\le S_r^{+}(t)$. This, together with the first part of the claim and the facts that $S_r^{+}(t)\ge S_r^{-}(t) = \Omega(n)$ (by \Cref{obs}(iv)) and $p_r\le p$ (by \Cref{obs}(ii)), implies the second part of the claim.
\end{proof}

\subsection{Upper bounding the maximum degree}
In this subsection, we show that for any fixed vertex $v$, w.v.h.p. \ref{upper bound degree} holds (see \Cref{claim:max-degree-glt} below). To do this, we need two intermediate results.

\begin{claim} \label{clm:bounded degree}
After Step~1 of the algorithm, w.v.h.p. any vertex $v\in V(G_r)$ has at most $\log^2 n$ neighbors in $\T^{\ell}$ for all $\ell \in [pS_r(0)]$.
\end{claim}

\begin{proof}
Define the random variable $X$ that counts the number of neighbors of $v$ in $\T^{\ell}$. Then $X$ is a sum of independent Bernoulli random variables. Note that
\[
\EE(X)=\sum_{i\in I^{\ell}_r(t)}\frac{\abs{N(v)\cap V^{\ell}_i(t)}}{\abs{V^{\ell}_i(t)}} \overset{\text{\ref{lower bound vijt}}}{\le} \sum_{i\in I^{\ell}_r(t)}\frac{\abs{N(v)\cap V^{\ell}_i(t)}}{S^{-}_r(t)} \overset{\text{\ref{upper bound degree}}}{\le} \frac{D_r(t)}{S^{-}_r(t)} \overset{\text{Obs. \ref{obs}(iii)}}{\le} 1.
\]
Hence, applying Chernoff's bound (the moreover part of \Cref{chernoff}) with $s=\log^2 n$, we obtain ${\PP[X \ge \log^2 n]\le \exp(-\log^2 n)}$, as desired.
\end{proof}

\begin{claim} \label{clm:bounded appearance}
After Step~1 of the algorithm, any vertex $v\in V(G_r)$ is contained in at most $\log^2 n$ sets $\T^{\ell}$ with $\ell \in [pS_r(0)]$.
\end{claim}

\begin{proof}
Let $X$ denote the number of $\ell \in [p S_r(0)]$ for which $v\in \T^{\ell}$. Then $X$ is a sum of independent Bernoulli random variables with mean
\[
\EE(X)=\sum_{\ell\in [pS_r(0)]}\PP[v\in \T^{\ell}] \le \sum_{\ell\in [pS_r(0)]}\frac{1}{\abs{V^{\ell}_i(t)}} \overset{\text{\ref{lower bound vijt}}}{\le} \frac{pS_r(0)}{S^{-}_r(t)}\overset{\text{Obs. \ref{obs}(iv)}}{=}O(p).
\]
Thus, applying Chernoff's bound (the moreover part of \Cref{chernoff}) with $s=\log^2 n$, we obtain ${\PP[X \ge \log^2 n]\le \exp(-\log^2 n)}$, as desired.
\end{proof}

We now use \Cref{clm:bounded degree,clm:bounded appearance} to show that each vertex loses the right proportion of its neighbors. This corresponds to the event \ref{upper bound degree} for an individual vertex $v$.

\begin{claim}\label{claim:max-degree-glt}
For any $\ell \in [pS_r(0)]$ and $v \in V(G_r)$ with $d_{G_r}(v,V^{\ell}(t))\ge \log^{15} n$, w.v.h.p. we have 
\[
d_{G_r}(v,V^{\ell}(t+1)) \leq (1-p+\e p/2)d_{G_r}(v,V^{\ell}(t)).
\]
In particular, for any $\ell \in [pS_r(0)]$ and $v \in V(G_r)$, w.v.h.p. we have $d_{G_r}(v,V^{\ell}(t+1)) \le D_r(t+1)$.
\end{claim}

\begin{proof}
Fix $\ell \in [pS_r(0)]$ and a vertex $v \in V(G_r)$ with $d_{G_r}(v,V^{\ell}(t))\ge \log^{15} n$. 
For $i\in I^{\ell}_r(t)$, let $c_i$ be the number of neighbors of $v$ in $V_i^{\ell}(t)$. Since the parts $V^{\ell}_i(t)$ with $i\in J^{\ell}$ were deleted in Step~4 of the algorithm, we find
\[
d_{G_r}(v,V^{\ell}(t+1))\le d_{G_r}(v,V^{\ell}(t))- \sum_{i\in I^{\ell}_r(t)} c_i \I_{i\in J^{\ell}}=\sum_{i\in I^{\ell}_r(t)} c_i \I_{i\notin J^{\ell}} =: Z.
\]
Here, we use $\I$ to denote the indicator random variable of the event given by the index. Thus, in order to prove the claim, it is enough to show that 
\begin{equation}\label{eq:Z}
 \PP\left[Z > (1-p+\e p/2)d_{G_r}(v,V^{\ell}(t))\right]\le \exp(-\omega(\log n)).   
\end{equation}
Recall that in Step~1 of the algorithm we randomly select vertices $\{v^s_i: s \in [pS_r(0)], i\in I^s_r(t)\}$, and then in Step~2 we activate each part $V^s_i(t)$ with probability $p$. As we have seen in \Cref{clm:bounded degree,clm:bounded appearance}, typically, the set of chosen vertices $\tilde{T}_s = \{v^s_i: i\in I^s_r(t)\}$ for each $s$ induces a graph with very small maximum degree, and no vertex lies in many of the sets $\tilde{T}_s$ with $s\in [pS_r(0)]$. Conditioning on this event when we activate the parts, we have a very small Lipschitz constant which allows us to use Talagrand's inequality.

Now, first expose the vertices $\{v^s_i: s \in [pS_r(0)], i\in I^s_r(t)\}$. Let $\mE$ be the intersection of the two events in 
\Cref{clm:bounded degree,clm:bounded appearance}. Note that $\mE$ is entirely determined by the choices of the $v^s_i$ in Step~1 where $s\in [pS_r(0)]$ and $i\in I_r^s(t)$, so all of the choices in Step~2 are still independent of it. By \Cref{clm:bounded degree,clm:bounded appearance}, we have

\begin{equation}\label{eq:E}
\PP[\overline{\mE}] \le \exp(-\omega(\log n)).    
\end{equation}
Next, we show that $Z$ is unlikely to be too large if $\mE$ holds. 
For this purpose, fix any choice $\T$ of $\{v^s_i\}$ for which $\mE$ holds. 
For $i\in I^{\ell}_r(t)$, let $Q_i$ be the set of $j\in I^{\ell}_r(t)$ with $v^{\ell}_i v^{\ell}_j \in E(G)$ and let $\bar{Q}_i$ be the set of $s \neq \ell$ for which $v^{\ell}_i\in \T^{s}$ (that is, $v^{\ell}_i=v^{s}_i$). From the assumption, we obtain $\abs{Q_i}\le \log^2 n$ and $\abs{\bar{Q}_i}\le \log^2 n$ for all $i\in I^{\ell}_r(t)$. Notice that $i\in J^{\ell}$ if and only if $i$ is activated for the $\ell$-th transversal but no $j\in Q_i$ is activated for the $\ell$-th transversal, and $i$ is not activated for any transversal $s\in \bar{Q}_i$. Thus, $\PP[i\in J^{\ell} \; | \; \T]\le p$.
Also, using $p=1/\log^3 n$, we have the following for every $i\in I^{\ell}_r(t)$. 
\[
\PP[i\in J^{\ell} \; | \; \T]\ge p\left(1-p|Q_i|-p|\bar{Q}_i|\right) \ge p(1-\eps/6).
\]
Thus, by the linearity of expectation, 
\begin{equation}\label{eq:left side of expected Z given T}
(1-p)d_{G_r}(v,V^{\ell}(t))\le \E(Z \; | \; \T)\le (1-p+\eps p/6)d_{G_r}(v,V^{\ell}(t)).
\end{equation}
We will use Talagrand's inequality to finish the proof. It is clear that $Z$ is $1$-certifiable since the event $i\notin J^{\ell}$ is witnessed by the non-activation of $i$ for the $\ell$-th transversal, or the activation of some $j\in Q_i$ for the $\ell$-th transversal, or the activation of $i$ for some transversal $s\in \bar{Q}_i$. 
To bound the Lipschitz constant of~$Z$, notice that whether $i\in I_r^{\ell}(t)$ is activated for the $\ell$-th transversal or not only affects the events $\{j\notin J^{\ell}\}$ with $j\in Q_i\cup \{i\}$. Also, for any $s\neq \ell$, whether $i \in I_r^s(t)$ is activated for the $s$-th transversal or not only affects the event $\{i\notin J^{\ell}\}$.  Since $G$ has local degree at most $C$, we have $0\le c_j \le C$ for all $j\in I_r^{\ell}(t)$. Therefore, for any $s\in [pS_r(0)]$, changing the decision whether $i\in I_r^s(t)$ is activated for the $s$-th transversal or not can affect $Z$ by an additive factor of at most $C(1+|Q_i|) \le 2C \log^2 n$. Thus, the Lipschitz constant of $Z$ is at most $2C \log^2 n$. 
We apply Talagrand's inequality with $s=\eps p\E(Z \; | \; \T)/6$ and $L=2C \log^2 n$. Noting that $60L\sqrt{\E(Z \; | \; \T)} \overset{\eqref{eq:left side of expected Z given T}}{\le} s \overset{\eqref{eq:left side of expected Z given T}}{\le} \eps pd_{G_r}(v,V^{\ell}(t))/6$ and $s^2/(L^2 \E(Z \; | \; \T))\overset{\eqref{eq:left side of expected Z given T}}{=}\Omega(\log^2 n)$ for $d_{G_r}(v,V^{\ell}(t))\ge \log^{15} n$, we obtain
\[
\PP[Z > (1-p+\eps p/2)d_{G_r}(v,V^{\ell}(t)) \; | \; \T] \le \exp(-\Omega(\log^2 n)). 
\]
As this bound holds for any choice $\T$ that satisfies $\mE$, we infer that
\[
\PP[Z > (1-p+\eps p/2)d_{G_r}(v,V^{\ell}(t)) \; | \; \mE] \le \exp(-\Omega(\log^2 n)).
\]
Together with \eqref{eq:E}, this implies \eqref{eq:Z}.
\end{proof}

\subsection{Many neighbors in partial transversals}

In this subsection, we show that for any individual vertex $v$, the events \ref{neighbor in transversals} and \ref{neighbor in remaining parts} hold w.v.h.p. (see \Cref{claim:neighbor in remaining parts,claim:neighbor in transversals} below). The following claim deals with  \ref{neighbor in remaining parts}. 

\begin{claim}\label{claim:neighbor in remaining parts}
For any vertex $v\in V(G_r)$, w.v.h.p. we have 
\[
\sum_{\ell \in [pS_r(0)]}\abs{N_{G_r}(v) \cap V^{\ell}(t+1)} \ge pS_r(0)\cdot (1- p-p_r-\delta p)^{t+1} d_{G_r}(v).
\]
\end{claim}

\noindent{\bf Remark.} The success of our randomized algorithm crucially rests on
the availability of the above claim. It is one of the main reasons for considering many transversals simultaneously.

\begin{proof}[Proof of \Cref{claim:neighbor in remaining parts}]
The claim is clearly true when $d_{G_r}(v)=0$. So suppose that $d_{G_r}(v)>0$.
Define the random variable $X:=\sum_{\ell}\abs{N_{G_r}(v) \cap V^{\ell}(t+1)}$. By the assumption \ref{neighbor in remaining parts}, after the $t$-th iteration we have
\[
\sum_{\ell}\abs{N_{G_r}(v) \cap V^{\ell}(t)} \ge pS_r(0)\cdot (1- p-p_r-\delta p)^t d_{G_r}(v).
\]
Consider any vertex $w\in \bigcup_{\ell} \left(N_{G_r}(v) \cap V^{\ell}(t)\right)$. 
Then $w\in N_{G_r}(v) \cap V^{\ell}_i(t)$ for some $\ell\in [p S_r(0)]$ and $i\in I^{\ell}_r(t)$. The entire part $V^{\ell}_i(t)$ gets deleted for the $\ell$-th transversal during the $(t+1)$-st iteration only if it was activated (with probability $p$) in Step~2, so this event occurs with probability at most $p$. 
The vertex $w$ gets deleted in Step~6 with probability $p_r$.
Thus, the probability that $w$ remains in $V^{\ell}(t+1)$ is at least $1-p-p_r$.
Hence, by the linearity of expectation, we have
\begin{equation}\label{eq:expectation of X}
\E(X) \ge (1-p-p_r)\sum_{\ell}\abs{N_{G_r}(v) \cap V^{\ell}(t)} \ge \left(1- p-p_r\right) \cdot pS_r(0)\cdot (1- p-p_r-\delta p)^t d_{G_r}(v).
\end{equation}
To prove concentration, for $\ell\in [pS_r(0)]$, define the random variable $X_{\ell} = \abs{N_{G_r}(v) \cap V^{\ell}(t+1)}$. Then $X=\sum_{\ell}X_{\ell}$. 
Note that the random variables $X_{\ell}$ are mutually independent since $X_{\ell}$ and also the set $N_{G_r}(v) \cap V^{\ell}(t+1)$ are completely determined by the random set $\hat{T}_{\ell}$ and the events $B_v^{\ell}$ for $v\in G_r$, and these are mutually independent for different $\ell$.
Moreover, $0\le X_{\ell} \le d_{G_r}(v)$. Thus, applying Hoeffding's inequality with $s=\delta p \E(X)$ and $L=d_{G_r}(v)$, and using \eqref{eq:expectation of X}, we obtain
\begin{align*}
\PP\left[X\le (1-\delta p)\E(X)\right]
& \le 2 \exp\left(-\delta^2 \left(1- p-p_r\right)^2p^3 S_r(0)\cdot (1- p-p_r-\delta p)^{2t} /2\right) \\
&= 2\exp(-\Omega(p^3n))\le \exp(-\Omega(\log^2 n)),
\end{align*}
where in the equality at the second line, we used the following estimates 
\[
(1- p-p_r)^2 \overset{\text{Obs. \ref{obs}}}{\ge} (1-2p)^2=\Omega(1), \; S_r(0)\overset{\text{Obs. \ref{obs}}}{=}\Omega(n), \; (1- p-p_r-\delta p)^{2t} \overset{\text{Obs. \ref{obs}}}{\ge}  (1-3p)^{2t} \overset{\text{Obs. \ref{obs}}}{=}\Omega(1).
\]
Hence, w.v.h.p. we have 
$X\ge (1-\delta p)\E(X) \overset{\eqref{eq:expectation of X}}\ge pS_r(0)\cdot (1- p-p_r-\delta p)^{t+1} d_{G_r}(v)$, as desired.
\end{proof}

In the remainder of this subsection, we deal with the event corresponding to \ref{neighbor in transversals}. We first show that each vertex has many neighbors in $\bigcup_{\ell}  \T^{\ell}$.

\begin{claim}\label{claim:neighbors in union of transversals}
After Step~1 of the algorithm, for any vertex $v\in V(G_r)$ with $d_{G_r}(v)\ge \log^5 n$, w.v.h.p. we have 
\[
\sum_{\ell \in [pS_r(0)]} \abs{N_{G_r}(v)\cap \T^{\ell}} \ge (1-\delta)p (1- (1+2\delta)p)^t d_{G_r}(v).
\]
\end{claim}

\begin{proof}
Let $X=\sum_{\ell} \abs{N_{G_r}(v)\cap \T^{\ell}}$. By the assumption \ref{neighbor in remaining parts}, after the $t$-th iteration we have
\[
\sum_{\ell}\abs{N_{G_r}(v) \cap V^{\ell}(t)} \ge pS_r(0)\cdot (1-p-p_r-\delta p)^t d_{G_r}(v).
\]
For any $\ell \in [p S_r(0)]$ and $i\in I^{\ell}_r(t)$, each vertex in $V^{\ell}_i(t)$ is included in $\T^{\ell}$ with probability $\frac{1}{\abs{V^{\ell}_i(t)}}$ which is, by \ref{lower bound vijt}, at least $\frac{1}{S_r^+(t)}$. Hence $\E(X) \ge \frac{1}{S_r^+(t)} \cdot \sum_{\ell}\abs{N_{G_r}(v) \cap V^{\ell}(t)}$. This inequality, combined with the above one, yields
\begin{align*}
\E(X)\ge \frac{pS_r(0)\cdot (1- p-p_r-\delta p)^t d_{G_r}(v)}{S_r^+(t)}&=p\left(\frac{1-p-p_r-\delta p}{1-p_r+p^2}\right)^td_{G_r}(v)\\
&\ge p (1- (1+2\delta)p)^t d_{G_r}(v)=\Omega(\log^2 n),
\end{align*}
where the last inequality holds since $\frac{1-p-p_r-\delta p}{1-p_r+p^2}=1-p-\frac{\delta+p_r+p(1-p)}{1-p_r+p^2}p \ge 1-p-2\delta p$ for $p,p_r=o(1)$ and $\delta=\Omega(1)$, and the last inequality follows from the facts that $p=1/\log^3 n$, $d_{G_r}(v) \ge \log^5 n$, and $(1-(1+2\delta)p)^t \ge (1-3p)^t\overset{\text{Obs. \ref{obs}(i)}}{=}\Omega(1)$.
Moreover, $X$ is a sum of independent Bernoulli random variables. Thus, applying Chernoff's bound with $s=\delta \E(X)$, and noting that $s^2/\E(X)=\delta^2\E(X)=\Omega(\log^2 n)$, we obtain

\begin{equation*}
\PP\left[X \leq (1-\delta) p (1- (1+2\delta)p)^t d_{G_r}(v)\right]\le \PP\left[X \leq (1-\delta) \E(X)\right]\le 2\exp\left(-\Omega(\log^2 n)\right). 
\end{equation*}
This proves \Cref{claim:neighbors in union of transversals}.
\end{proof}

We are now ready to take on \ref{neighbor in transversals}.

\begin{claim}\label{claim:neighbor in transversals}
For any vertex $v \in V(G_r)$ with $d_{G_r}(v)\ge \log^{15} n$, w.v.h.p. we have
\[
\abs{N_{G_r}(v) \cap T(t+1)} \ge  (1-3\delta) p^2\sum_{0\le j\le t} (1-(1+2\delta)p)^{j} d_{G_r}(v).
\]
\end{claim}

\begin{proof}
From the definition of $T(t+1)$, we see that
\[
\abs{N_{G_r}(v) \cap T(t+1)}=\abs{N_{G_r}(v) \cap T(t)}+\sum_{\ell\in [pS_r(0)]} \abs{N_{G_r}(v) \cap \bT^{\ell}}.  
\]
Moreover, by the assumption \ref{neighbor in transversals},
$\abs{N_{G_r}(v) \cap T(t)} \ge  (1-3\delta)p^2 \sum_{0\le j\le t-1} (1-(1+2\delta)p)^j d_{G_r}(v)$.
Thus, to prove the claim, it suffices to show that w.v.h.p.
\[
Z:= \sum_{\ell\in [pS_r(0)]} \abs{N_{G_r}(v) \cap \bT^{\ell}} \ge  (1-3\delta)p^2 (1-(1+2\delta)p)^t d_{G_r}(v).
\]
Similar to the proof of \Cref{claim:max-degree-glt}, we first randomly choose the vertices $\{v^\ell_i: \ell \in [pS_r(0)], i\in I^\ell_r(t)\}$ in Step~1, and let $\mE$ be the intersection of the events in 
Claims~\ref{clm:bounded degree}, \ref{clm:bounded appearance}, and~\ref{claim:neighbors in union of transversals}.
Then $\PP[\overline{\mE}] \le \exp(-\omega(\log n))$.    
Now, fix any choice $\T$ of $\{v^\ell_i\}$ for which $\mE$ holds. By \Cref{claim:neighbors in union of transversals}, 
\[
\sum_{\ell \in [pS_r(0)]} \abs{N_{G_r}(v)\cap \T^{\ell}} \ge (1-\delta)p (1- (1+2\delta)p)^t d_{G_r}(v).
\]
Note that 
\begin{equation}\label{eq:double summation for Z}
Z = \sum_{\ell\in [pS_r(0)]} \sum_{w\in N_{G_r}(v)\cap \T^{\ell}} \I_{w\in \bT^{\ell}}.
\end{equation}
To this end, for every $w\in N_{G_r}(v)\cap \T^{\ell}$, we will estimate the probability that $w\in \bT^{\ell}$. Notice that for any vertex $w\in N_{G_r}(v)\cap \T^{\ell}$, we must have $w=v^{\ell}_i$ for some $i \in I^{\ell}_r(t)$. 
For $\ell\in [pS_r(0)]$ and $i\in I^{\ell}_r(t)$, let $Q_i^{\ell}$ be the set of $j\in I^{\ell}_r(t)$ with $v^{\ell}_i v^{\ell}_j \in E(G)$ and let $\bar{Q}_i^{\ell}$ be the set of $s \neq \ell$ for which $v^{\ell}_i\in \T^{s}$ (that is, $v^{\ell}_i=v^{s}_i$). From the assumption, we obtain $\abs{Q_i^{\ell}}\le \log^2 n$ and $\abs{\bar{Q}_i^{\ell}}\le \log^2 n$ for all $\ell\in [pS_r(0)]$ and $i\in I^{\ell}_r(t)$. Let $v^{\ell}_i\in N_{G_r}(v)\cap \T^{\ell}$ for some $\ell \in [pS_r(0)]$ and $i \in I^{\ell}_r(t)$. Then, notice that $v^{\ell}_i\in \bT^{\ell}$ if and only if $i\in J^{\ell}$ if and only if $i$ is activated for the $\ell$-th transversal but no $j\in Q_i^{\ell}$ is activated for the $\ell$-th transversal, and $i$ is not activated for any transversal $s\in \bar{Q}_i^{\ell}$. 
Thus, as in the proof of \Cref{claim:max-degree-glt}, for every $\ell\in [pS_r(0)]$ and $i\in I^{\ell}_r(t)$, we have $\PP[i\in J^{\ell} \; | \; \T]\ge p\left(1-p|Q_i^{\ell}|-p|\bar{Q}_i^{\ell}|\right) \ge (1-\delta)p$, and so the probability that $v^{\ell}_i$ remains in $\bT^{\ell}$ is at least $(1-\delta)p$.
Thus, using \eqref{eq:double summation for Z} and the linearity of expectation, we obtain
\begin{align}\label{eq:Z given T}
\E[Z \; | \; \T]\ge (1-\delta)p\cdot\sum_{\ell \in [pS_r(0)]} \abs{N_{G_r}(v)\cap \T^{\ell}} &\ge (1-\delta)p \cdot (1-\delta)p(1- (1+2\delta)p)^t d_{G_r}(v) \nonumber \\ 
&\ge (1-2\delta)p^2 (1- (1+2\delta)p)^t d_{G_r}(v).
\end{align}
It is clear that $Z$ is $(2\log^2 n +1)$-certifiable since the event $v^{\ell}_i\in \bT^{\ell}$ (equivalently, $i\in J^{\ell}$) is witnessed by the activation of $i$ for the $\ell$-th transversal, the non-activation of $j$ for every $j\in Q_i^{\ell}$ for the $\ell$-th transversal, and the non-activation of $i$ for every transversal $s\in \bar{Q}_i^{\ell}$. 
To bound the Lipschitz constant of~$Z$, notice that whether $i\in I_r^{\ell}(t)$ is activated for the $\ell$-th transversal or not only affects the events $\{j\in J^{\ell}\}$ with $j\in Q_i^{\ell}\cup \{i\}$, and the events $\{i\in J^{s}\}$ with $s\in \bar{Q}_i^{\ell}$. 
Therefore, for any $\ell\in [pS_r(0)]$ and $i\in I_r^{\ell}(t)$, changing the decision whether $i$ is activated for the $\ell$-th transversal or not can affect $Z$ by an additive factor of at most $1+|Q_i^{\ell}|+|\bar{Q}_i^{\ell}| \le 3\log^2 n$. Thus, the Lipschitz constant of $Z$ is at most $3\log^2 n$. 
Applying Talagrand's inequality with $s=(\delta/2)\E(Z \; | \; \T)$ and $L=r=3\log^2 n$, and noting that $60L\sqrt{r\E(Z \; | \; \T)} \le s$ and $s^2/(L^2 r\E(Z \; | \; \T))=\Omega(\log^2 n)$ (to get these, we use \eqref{eq:Z given T}, \Cref{obs}(i), and the fact that $d_{G_r}(v)\ge \log^{15} n$), we obtain
\[
\PP[Z < (1-\delta)\E(Z \; | \; \T) \; | \; \T] \le \exp(-\Omega(\log^2 n)). 
\]
As this bound holds for any choice $\T$ that satisfies $\mE$, using \eqref{eq:Z given T} and the fact that $\mE$ holds w.v.h.p., it follows that w.v.h.p. 
\begin{equation*}
Z \ge (1-\delta)\cdot (1-2\delta)p^2 (1- (1+2\delta)p)^t d_{G_r}(v) \ge (1-3\delta)p^2 (1- (1+2\delta)p)^t d_{G_r}(v). \qedhere
\end{equation*}
\end{proof}

\subsection{Wrapping up the proof}
We now use Claims~\ref{claim:B4}, \ref{claim:max-degree-glt}, \ref{claim:neighbor in remaining parts}, and~\ref{claim:neighbor in transversals} to show that with positive probability the events \ref{lower bound vijt}--\ref{neighbor in remaining parts} {\em simultaneously} hold, thereby completing the proof of \Cref{thm:main}.  We cannot simply use a union bound since we do not assume any bound on the number of parts $k$.

\begin{lemma} \label{lem:apply local lemma}
With positive probability, after the $(t+1)$-st iteration of the $r$-th round, the events \ref{lower bound vijt}--\ref{neighbor in remaining parts} for $t+1$ simultaneously hold.
\end{lemma}

\begin{proof}
To prove the assertion, we use the Lov\'asz Local Lemma.
For each $i\in [k]$, consider the following four events.
\begin{itemize}[leftmargin=*]
    \item For $i\in I_r^{\ell}(t+1)$, we have $S_r^{-}(t+1)\le \abs{V_i^{\ell}(t+1)}\le S_r^{+}(t+1)$. 
    \item For all $v\in V_i\cap V(G_r)$ and $\ell\in [p S_r(0)]$, there are at most $D_r(t+1)$ edges from $v$ to $V^{\ell}(t+1)$.
    \item For every $v\in V_i\cap V(G_r)$ with $d_{G_r}(v)\ge \log^{15} n$, we have 
    \[
    \abs{N_{G_r}(v) \cap T(t+1)} \ge  (1-3\delta) p^2 \sum_{0\le j< t+1} (1-(1+2\delta)p)^j d_{G_r}(v).
    \]
    \item For every $v\in V_i\cap V(G_r)$, we have 
    \[
    \sum_{\ell\in [p S_r(0)]}\abs{N_{G_r}(v) \cap V^{\ell}(t+1)} \ge pS_r(0)\cdot (1- p-p_r-\delta p)^{t+1} d_{G_r}(v).
    \]
\end{itemize}
Let $\B_i$ denote the union of the complements of these events. Since for any $i\in [k]$, there are at most $n^2$ pairs $(v,\ell)$ where $v\in V_i$ and $\ell\in [pS_r(0)]$, Claims~\ref{claim:B4}, \ref{claim:max-degree-glt}, \ref{claim:neighbor in remaining parts}, and~\ref{claim:neighbor in transversals} together with the union bound imply that $\PP[\B_i]=\exp(-\omega(\log n))$. Now consider a graph $\G$ on the vertex set $[k]$ where $x,y\in [k]$ are joined if there is an edge (in $G$) between $V_x$ and $V_{y}$. Observe that each event $\B_i$ is fully determined by the random choices (including the Bernoulli random variables $B^{\ell}_v$) involving the parts $V_{j}$ with $\dist_{\G}(i,j)\le 2$.  
Thus, $\G^4$, the graph obtained from $\G$ by joining pairs of vertices of distance at most $4$, is a dependency graph for the events $\{\B_i\}_{i\in [k]}$. Moreover, since $\Delta(G)\le n$ and $|V_i|=n$ for all $i\in [k]$, we have $\Delta(\G)=O(n^2)$, implying $\Delta(\G^4)=O(n^8)$. Hence, the condition of the Lov\'asz Local Lemma is fulfilled and so with positive probability none of the events $\B_i$ occurs. In other words, there is an outcome that satisfies the events \ref{lower bound vijt}--\ref{neighbor in remaining parts} for $t+1$.
\end{proof}
This finishes the proof of \Cref{thm:main}.

\section{Reducing local degree} \label{sec:reduction}
In this section, we prove \Cref{thm} assuming \Cref{thm:main}. We make use of the following two lemmas, which are generalizations of Lemmas~3.2 and~3.3 in \cite{LS07}. Although their proofs can be easily modified to obtain these generalizations, we include the details for the sake of completeness.

\begin{lemma} \label{lem:nice partition}
For any $0<\e<1$, there exists $n_0$ such that the following holds for all $n \ge n_0$. Let $G$ be a multipartite graph with maximum degree at most $(1-\e)n$, parts $V_1,\ldots,V_k$ of size $|V_i| = n$, and local degree at most $n^{1/3}$. Let $m=n^{2/3}$. Then, there exist partitions $V_i=\bigcup_{\ell\in [m]} V_i^{\ell}$ for $i \in [k]$ with the following properties.
\stepcounter{propcounter}
\begin{enumerate}[leftmargin=*,label = {\bfseries \emph{\Alph{propcounter}\arabic{enumi}}}]
\item For every $i\in [k]$ and every $\ell\in [m]$, we have $\abs{V_i^{\ell}}\ge (1-\e/4)n^{1/3}$. 
\label{vij size}
\item For every $\ell\in [m]$, each vertex $v$ has at most $(1-3\e/4)n^{1/3}$ neighbors in $\bigcup_{i\in [k]} V_i^{\ell}$. 
\label{vij degree}
\item For every $i\in [k]$ and every $\ell\in [m]$, each vertex $v$ has less than 12 neighbors in $V_i^{\ell}$. \label{vij local degree}
\end{enumerate}
\end{lemma}

\begin{proof}
For each $i\in [k]$, we form a partition $V_i=\bigcup_{\ell\in [m]}V_i^{\ell}$ of $V_i$ by randomly assigning each $v\in V_i$ into one of the $m$ sets with equal probability $p=1/m$. 
We define the following three types of bad events corresponding to each of \ref{vij size}--\ref{vij local degree}. 
For each $i\in [k]$, let $\mA_i$ be event that $\abs{V_i^{\ell}}<(1-\e/4)n^{1/3}$ for some $\ell\in [m]$. 
For each vertex $v$, let $\mB_v$ be the event that for some $\ell\in [m]$, the number of neighbors of $v$ in $\bigcup_{i\in [k]} V_i^{\ell}$ exceeds $(1-3\e/4)n^{1/3}$.
Finally, for each vertex $v$ and part $V_i$, let $\mC_{v,i}$ be the event that for some $\ell\in [m]$, the number of neighbors of $v$ in $V_i^{\ell}$ is at least $12$. 
We use the local lemma to show that with positive probability, none of these bad events occur. 

Each of the events $\mB_v,\mC_{v,i}$ is completely determined by the choices for neighbors of $v$, and $\mA_i$ is completely determined by the choices for vertices in $V_i$. Since the maximum degree of $G$ is at most $n$ and $|V_i|=n$ for all $i$, each event is mutually independent of all but $O(n^2)$ other events.

We next upper bound the probabilities of bad events. Since the size of $V_i^{\ell}$ is binomially distributed with mean $|V_i|p=n^{1/3}$, it follows from Chernoff's bound that the probability that $\abs{V_i^{\ell}}<(1-\e/4)n^{1/3}$ is at most $\exp(-\Omega(n^{1/3}))$. Hence, by the union bound, $\PP[\mA_i] \le m\cdot \exp(-\Omega(n^{1/3})) \le n^{-3}$.
Similarly, since the number of neighbors of a vertex $v$ in $\bigcup_{i\in [r]} V_i^{\ell}$ is binomially distributed with mean at most
$(1-\e)n\cdot p=(1-\eps)n^{1/3}$, Chernoff's bound implies that the probability that it exceeds $(1-3\e/4)n^{1/3}$ is at most $\exp(-\Omega(n^{1/3}))$. Thus, using the union bound again, we find $\PP[\mB_v] \le m\cdot \exp(-\Omega(n^{1/3})) \le n^{-3}$. Finally, as $v$ has at most $n^{1/3}$ neighbors in $V_i$, we have
\[
\PP[\mC_{v,i}] \le m\cdot \binom{n^{1/3}}{12}p^{12} \le n^{-3}.
\]
Therefore, by the Lov\'asz Local Lemma, with positive probability none of the events $\mA_i,\mB_{v},\mC_{v,i}$ occur, and we obtain partitions satisfying conditions \ref{vij size}--\ref{vij local degree}.
\end{proof}

The following lemma simultaneously generalizes Lemma~3.3 in~\cite{A92} and Lemma~3.3 in~\cite{LS07}.  

\begin{lemma} \label{lem: nice bipartition}
The following holds for $\D$ sufficiently large and $d>\log^4 \D$. Let $G$ be a multipartite graph with maximum degree at most $\D$, parts $V_1,\ldots,V_k$ of size $|V_i| = 2n$, and local degree at most $d$. Then, there exist partitions $V_i=V_i^1\cup V_i^2$ for $i \in [k]$ such that 
\stepcounter{propcounter}
\begin{enumerate}[leftmargin=*,label = {\bfseries \emph{\Alph{propcounter}\arabic{enumi}}}]
\item For every $i\in [k]$, we have $\abs{V_i^1}=\abs{V_i^2}=n$.
\item For every $\ell\in [2]$, every vertex $v\in \bigcup_{i\in [k]} V_i^{\ell}$ has at most $\D/2 + \D^{2/3}$ neighbors in $\bigcup_{i\in [k]} V_i^{\ell}$. \label{vi2 degree}
\item For every $i\in [k]$ and every $\ell\in [2]$, every vertex $v\in \bigcup_{j\in [k]} V_j^{\ell}$ has at most $d/2 + d^{2/3}$ neighbors in $V_i^{\ell}$. \label{vi2 local degree}
\end{enumerate}
\end{lemma}

\begin{proof}
Arbitrarily pair up the vertices in each $V_i$. For each pair of vertices in $V_i$, designate one of them randomly and independently of other events to $V_i^1$ and the other one to $V_i^2$. By construction, $\abs{V_i^1}=\abs{V_i^2}=n$ for all $i$. 
For each vertex $v$, define $\mA_v$ to be the event that for some $\ell\in [2]$, the number of neighbors of $v$ in $\bigcup_{i\in [k]}V_i^{\ell}$ exceeds $\D/2+\D^{2/3}$. For each vertex $v$ and part $V_i$, let $\mB_{v,i}$ be the event that for some $\ell\in [2]$, the number of neighbors of $v$ in $V_i^{\ell}$ exceeds $d/2+d^{2/3}$. We use the local lemma to show that none of the bad events $\mA_v, \mB_{v,i}$ happen.

By the same argument as in the proof of \Cref{lem:nice partition}, we can show that each bad event is mutually independent of all but at most $O(\D^2)$ such events. We next bound the probabilities of bad events.
Consider event $\mA_v$. Observe that if two neighbors of $v$ are paired up, then exactly one of them will lie in $\bigcup_{i\in [k]} V_i^{\ell}$. Let $T$ be the set of all neighbors of $v$ that are paired to vertices that are not neighbors of $v$. Let $x$ denote the number of vertices in $T$ that are designated to $\bigcup_{i\in [k]} V_i^{\ell}$. 
Then the number of neighbors of $v$ in $\bigcup_{i\in [k]} V_i^{\ell}$ is upper bounded by $x + (\D-|T|)/2$. The quantity $x$ is binomially distributed with parameters $|T| \le \D$ and $1/2$. By Chernoff's bound, the probability that it exceeds $|T|/2+\D^{2/3}$ is bounded by $2 \exp(-2(\D^{2/3})^2/(3|T|)) \le \D^{-3}$. Hence $\PP[\mA_v] \le \D^{-3}$. Similarly, one can use the assumption $d>\log^4\D$ to show $\PP[\mB_{v,i}]\le \D^{-3}$. Therefore, by the Lov\'asz Local Lemma, with positive probability, none of the events $\mA_v,\mB_{v,i}$ occur. 
\end{proof}

Using the previous lemmas and \Cref{thm:main}, we finally prove our main result \Cref{thm}. The argument is again similar to \cite{LS07}.

\begin{proof}[Proof of \Cref{thm}]
Let $0< \eps < 1$ and $\gamma>0$ be a sufficiently small constant dependent on $\eps$. Let $G$ be a multipartite graph with maximum degree at most $(1- \eps)n$, parts $V_1,\dots,V_k$ of size $|V_i|\ge n$, and local degree at most $\gamma n$. Since deleting vertices from $G$ does not increase the maximum degree or the local degree, we may assume that $|V_i|=n$ for all $i\in [k]$. Since there is nothing to prove when $n < \gamma^{-1}$, we assume that $n\ge \gamma^{-1}$. Let $m = n^{2/3}$. 

First, we consider the case $n\le \gamma^{-4/3}$. Then, the local degree is at most $\gamma n \le n^{1/4}$. 
Thus, applying \Cref{lem:nice partition} gives us partitions $V_i = \bigcup_{\ell\in [m]} V_i^{\ell}$ for $i\in [k]$ satisfying \ref{vij size}--\ref{vij local degree}. 
Therefore, for every $\ell \in [m]$, we apply \Cref{thm:main} on the multipartite graph induced by $V_i^{\ell}, \dots, V_k^{\ell}$ with $n$ replaced by $(1- \eps/4)n^{1/3}$ and $\eps$ replaced by $\eps/2$ to get $(1-\eps)n^{1/3}$ pairwise disjoint independent transversals. (This application is possible since $(1-\eps/4) n^{1/3}\ge (1-\eps/4) \gamma^{-1/3}$ and $\gamma$ is sufficiently small.) 
Taking the union of these independent transversals over $\ell \in [m]$ gives us $(1-\eps)n$ pairwise disjoint independent transversals in $G$, as desired. 

We now assume $n > \gamma^{-4/3}$. Let $j\ge 1$ be the integer such that $2^{j-1} < \gamma^{4/3} n \le 2^j$. We now delete at most $2^j$ vertices from each $V_i$ to ensure that $2^j$ divides $|V_i|$ for $i\in [k]$. Then clearly, $|V_i|\ge n- 2^j$ for $i\in [k]$. 
Define $\D_0 = (1-\eps)n$ and $d_0 = \gamma n$, and for every $t\ge 0$, define 
\[
\D_{t+1} = \frac{\D_t}{2} + \D_t^{2/3} \;\;\;\; \text{and} \;\;\;\; d_{t+1} = \frac{d_t}{2} + d_t^{2/3}.
\] 
We claim the following. 
\stepcounter{propcounter}
\begin{enumerate}[leftmargin=*,label = {\bfseries \emph{\Alph{propcounter}\arabic{enumi}}}]
\item $\frac{1}{4\gamma^{4/3}} < \D_j \le (1-\eps/2)\left(\frac{n}{2^j} -1\right)$, \label{max degree bound}
\item $d_j \le \left(\frac{n}{2^j} - 1\right)^{1/3}$, and \label{local degree upper bound}
\item $d_t > \log^4 \D_t$ for $0\le t < j$. \label{local degree not too small}
\end{enumerate}
We first finish the proof using these claims and prove them at the end. Let $n' = |V_i|/2^j$ for $i\in [k]$. 
Then, $\frac{n}{2^j} - 1\le n'\le \frac{n}{2^j}$.
Using \ref{local degree not too small} and applying \Cref{lem: nice bipartition} iteratively $2^j -1$ times, we get a partition of $V_i$ into $V_i^1,\dots,V_i^{2^j}$ such that for every $i\in [k]$ and $\ell\in [2^j]$, we have $\abs{V_i^{\ell}} = |V_i|/2^j = n'$ and \ref{vi2 degree}, \ref{vi2 local degree} hold with $\D$ and $d$ replaced by $\D_j$ and $d_j$. 
Moreover, since \ref{max degree bound} and \ref{local degree upper bound} hold, we have $\D_j\le (1-\eps/2)n'$ and $d_j\le n'^{1/3}$.
Therefore, we apply \Cref{lem:nice partition} with $n$ replaced by $n'$ and $\eps$ replaced by $\eps/2$ to get partitions $V_i^{\ell} = \bigcup_{\ell'\in [n'^{2/3}]} V_i^{\ell, \ell'}$ for $i\in [k], \ell\in [2^j]$ satisfying \ref{vij size}--\ref{vij local degree} with $n$ replaced by $n'$ and $\eps$ replaced by $\eps/2$. 
Next, as in the case $n\le \gamma^{-4/3}$, apply \Cref{thm:main} to each multipartite graph induced by $V_1^{\ell,\ell'} \cup \dots \cup V_k^{\ell,\ell'}$ to get $(1-\eps)n$ pairwise disjoint independent transversals in $G$.

We now return to proving \ref{max degree bound}--\ref{local degree not too small}. 
For the lower bound in \ref{max degree bound}, notice that for every $0\le t\le j$,
\begin{equation}\label{eq:lower bound on Delta t}
\D_t \ge \frac{\D_0}{2^t} = \frac{(1-\eps)n}{2^t} \ge \frac{(1-\eps)n}{2^j} > \frac{1-\eps}{2\gamma^{4/3}} > \frac{1}{4\gamma^{4/3}}.
\end{equation}
To prove the upper bound in \ref{max degree bound}, first note that 
\[
\D_{t+1} = \frac{\D_t}{2} + \D_t^{2/3}\le \frac{1}{2}(\D_t^{1/3} + 1)^3 \;\;\;\; \text{and thus} \;\;\;\; \D_{t+1}^{1/3} \le \frac{\D_t^{1/3}}{2^{1/3}} + \frac{1}{2^{1/3}}.
\]
Consequently, 
\[
\D_{j}^{1/3} \le \frac{\D_0^{1/3}}{2^{j/3}} + \sum_{t\in [j]} \frac{1}{2^{t/3}} \le \frac{((1-\eps)n)^{1/3}}{2^{j/3}} + 4 \le (1+\eps/4)^{1/3}\frac{((1-\eps)n)^{1/3}}{2^{j/3}},
\]
where the last inequality uses $n/2^j > \gamma^{-4/3}/2$ and that $\gamma$ is sufficiently small. 
Thus,
\[
\D_j \le \frac{(1+\eps/4)(1-\eps)n}{2^j} \le (1-\eps/2)\left(\frac{n}{2^j} -1\right).
\]
It follows that for every $0\le t < j$, we have 
\begin{equation} \label{eq:bound on Delta}
\D_t \le 2\D_{t+1}\le \cdots\le 2^{j-t}\D_j \le (1-\eps/2) \frac{n}{2^t}.
\end{equation}

To prove \ref{local degree upper bound}, by dealing with $d_t$ the same way we dealt with $\D_t$, we obtain 
\[
d_{j}^{1/3} \le \frac{(\gamma n)^{1/3}}{2^{j/3}} + 4 \le 2 \frac{(\gamma n)^{1/3}}{2^{j/3}},
\]
where the last inequality uses $\gamma n/2^j > \gamma^{-1/3}/2$ and that $\gamma$ is sufficiently small. By the choice of $j$, we have $\left(\frac{2^{j-1}}{n}\right)^{3/4}< \gamma \le \left(\frac{2^j}{n}\right)^{3/4}$ and thus, since $\gamma$ is chosen sufficiently small, $n/2^j$ is sufficiently large. Using these facts, we have 
\[
d_j \le \frac{8\gamma n}{2^j} \le 8 \left(\frac{n}{2^j}\right)^{1/4} \le \left(\frac{n}{2^j} - 1\right)^{1/3}.
\]

Finally, to prove \ref{local degree not too small}, note that for every $0\le t < j$, we have $\gamma \ge \left(\frac{2^t}{n}\right)^{3/4}$. Thus, using \eqref{eq:bound on Delta}, we see that
\begin{equation*}
d_t \ge \frac{\gamma n}{2^t} \ge \left(\frac{n}{2^t}\right)^{1/4} \ge \left(\frac{\D_t}{1-\eps/2}\right)^{1/4} \ge \log^4 \D_t, 
\end{equation*}
where the last inequality uses the fact that $\Delta_t$ is sufficiently large, which easily follows from \eqref{eq:lower bound on Delta t} and the fact that $\gamma$ is sufficiently small.
\end{proof}

\section{Concluding remarks} \label{sec:remarks}
Now that we have an asymptotic solution to Loh--Sudakov's conjecture in \Cref{thm}, it seems promising to try and turn it into an exact solution using absorption.
There are absorption lemmas for similar problems; see, e.g., \cite{LM13}. However, in these results, the number of parts is fixed and does not seem to extend to our setting.

It was suggested in \cite{GS22} to develop a hypergraph version of \Cref{Loh} or its extension in \cite{GS22}. In a similar vein, it would be interesting to find a hypergraph generalization of our main result \Cref{thm}.

\section*{Acknowledgements}
We thank the anonymous referee for their helpful comments, as well as Tianjiao Dai and Xin Zhang for pointing out a typo in the statement of Talagrand's inequality used in a previous version of this paper.

\bibliographystyle{abbrv}
\bibliography{references.bib}

\end{document}